\documentclass[dvips,epsfig,reqno,graphics]{amsart}
\usepackage{amssymb}
\usepackage{amsfonts}
\usepackage{amsmath}
\usepackage{amsxtra}

\title[The geometry of the Second Fundamental Form]
{On the Geometry of the Second Fundamental Form of Translation Surfaces in $\mathbb{E}^3$}
\author{Marian Ioan Munteanu}
\author{Ana Irina Nistor}
\thanks{The first author was supported by
grant PN-II ID 398/2007-2010, CNCSIS, Romania. The second author was partially supported by
grant PN-II ID xx/2008-2011, CNCSIS, Romania.}

\newtheorem{contor}{1.}
\newtheorem{proposition}[contor]{Proposition}
\newtheorem{theorem}[contor]{Theorem}

\def\proof{{\sc Proof.\ }}

\def\R{\mathbb{R}}
\def\E{\mathbb{E}}

\date{ \today }

\begin{document}

\maketitle

\input epsfx.tex

\begin{abstract}
In this paper we study the second fundamental form of translation surfaces in
$\mathbb{E}^3$. We give a non-existence result for polynomial translation surfaces in $\E^3$
with vanishing second Gaussian curvature $K_{II}$.
We classify those translation surfaces for which $K_{II}$ and $H$ are proportional.
Finally we obtain that there are no $II$--minimal translation surfaces in the Euclidean 3-space.

\vspace{1mm}


\vspace{1mm}


\end{abstract}

\section{ Introduction}

Given a surface $M$ immersed in Euclidean 3-space, the knowledge of its first
fundamental form $I$ and its second fundamental form $II$
facilitates the analysis and the classification of surface shape.
$I$ is an intrinsic object of the surface that measures the amount of movement
of $M$ at a point of the surface in which is expressed, being invariant
to translations and rotations of the surface in ambient 3-space. Concerning the
second fundamental form $II$, it is an extrinsic tool to characterize the twist of
the surface in the ambient. However, it is dependent on the position of the surface
in the ambient $3$-space.

\smallskip
Over the years, the geometry of the second fundamental form $II$ has been a
popular topic for many  researchers. Very recent results concerning the
curvature properties associated to $II$ and other variational aspects
can be found in \cite{kn:Ver08}. It is a basic fact that we can talk about
characteristics of $M$ measured by means of its second fundamental form only
when it can be considered a metric tensor on the surface.
It is easy to see that $II$ is a metric tensor on $M$ if and only if it is non-degenerate.
From \cite{kn:Sod05} we get the following criterion for
non-degeneracy of $II$: {\em The second fundamental form $II$ of M is non-degenerate if and
only if $M$ is non-developable.}

\smallskip
At this point, on a non-developable surface $M$ we can consider the second fundamental
form $II$ as a new Riemannian metric on $M$ and the second
Gaussian curvature (respectively second mean curvature), denoted by $K_{II}$ (respectively by $H_{II}$),
is nothing else but the Gaussian curvature
(resp., mean curvature) of $(M, II)$. Several formulae for $K_{II}$, or in generally
{\em the scalar curvature of the second fundamental form} of a
hypersurface in a Semi-Riemannian space with a non-degenerate $II$, in different
ambient spaces can de found in \cite{kn:Ver08}. (See also References in
\cite{kn:Ver08}.) Regarding the second mean curvature, the critical points of the
area functional of the second fundamental form are those surfaces for
which the mean curvature of the second fundamental form $H_{II}$ vanishes (see \cite{kn:Ver08}).

An extended study was made by several authors on the second Gaussian curvature.
For example, in the case of the surfaces of revolution D. Koutroufiotis
has shown in \cite{kn:Kou74} that a closed ovaloid is a sphere if
$K_{II}=cK,\ c\in\mathbb{R}$ or if $K_{II}=\sqrt{K}$. Another property of spheres
was proved by T. Koufogiorgos and T. Hasanis in \cite{kn:KH77} meaning
that the sphere is the only closed ovaloid satisfying $K_{II}=H$.  Since
a natural generalization of these surfaces are the helicoidal surfaces,
the property $K_{II}=H$ was studied also for this surfaces. Therefore,
C. Baikoussis and T. Koufogiorgos in \cite{kn:BK97} proved that the helicoidal
surfaces satisfying $K_{II}=H$ are locally characterized by constant
ratio of the principal curvatures. Concerning the class of non-developable ruled surfaces,
one first result was obtained by D.E. Blair and
Th. Koufogiorgos in \cite{kn:BK92} studying relations between $K_{II}$ and $H$ for this
surfaces in $\mathbb{E}^3$. Then, the study was extended by
Y.H. Kim and D.W. Yoon in \cite{kn:KY04}, W. Sodsiri in \cite{kn:Sod05}, D.W. Yoon in
\cite{kn:Yoo06} for 3--dimensional Lorentz--Minkowski spaces and
for different relations between $H$, $K$, $H_{II}$ and $K_{II}$.

Previous results in the case of translation surfaces, the subject of our paper,
can be found in \cite{kn:GW07}. The authors give the characterization
theorem for translation surfaces with vanishing second Gaussian curvature in
$\mathbb{E}^3$ and $\mathbb{E}_1^{3}$ as follows:

{\em {\bf Theorem.}
If a translation surface in $\mathbb{E}_1^{3}$ parametrized by
$
\bar{x}(u,v)=(u,v, f(u)+g(v))
$ has $K_{II}=0$, then
\begin{center}
$
f(u)=\int F^{-1}(u+d)du\ \ $
and $\  \ g(v)=\int G^{-1}(v+m)dv$
\end{center}
with $F$ and $G$ real functions determined by
\begin{center}
$F(x)=\int \frac{x^2}{ax^4+bx^2+c}dx$
and
$G(x)= \int \frac{x^2}{-ax^4+(2a+b)x^2-a-b-c} dx$,
\end{center}

and $a,b,c,d$ \c si $m$ real numbers.}

\medskip
With all these results in mind, two important properties follow immediately as it is shown
in the definition mentioned in \cite{kn:Sod05}:
{\em A non-developable surface is said to be $II-flat$ if $K_{II}=0$ and respectively $II-minimal$ if $H_{II}=0$}.

The aim of this article is to study the {\em translation surfaces}
in Euclidean 3-space from the point of view of their second fundamental form, namely
their properties involving $K_{II}$ and $H_{II}$. Also some important
observations concerning a particular type of translation surfaces,
the {\em polynomial translation surfaces} in $\E^3$ will be made.
By definition, a translation surface in $\mathbb{E}^3$ is given locally by an immersion
$r:U\subseteq \mathbb{R}^2\rightarrow \mathbb{R}^3$ such that
$(u,v)\mapsto(u,v,f(u)+g(v))$, where $f$ and $g$ are smooth functions.
These surfaces are important either because they are
interesting themselves or because they furnish counterexamples for some problems
(e.g. it is a known fact that a minimal surface has
vanishing second Gaussian curvature but not conversely -- see for details \cite{kn:BK92}).

\smallskip
For the rest of this paper we call {\em polynomial translation surfaces}
(in short, {\em PT surfaces}) those translation surfaces
for which $f$ and $g$ are polynomials (see also \cite{kn:MN08}).
For technological applications in which different surfaces are needed (such as
Computed Aided Manufacturing) polynomial forms are preferred since they may be
incorporated into the CAD software in order to be easily
processed by numerical computations.

\smallskip
Next section is dedicated to the study of the {\em II-flat PT surfaces},
where we formulate the classification theorem. Section 3 includes the main
theorems concerning the translation surfaces and PT surfaces which satisfy a
relation of type $K_{II}=\lambda H$. These surfaces are called
{\em generalized Weingarten surfaces} (See for details \cite{kn:Sod05}).
The last section presents some results involving the
II-minimality property of translation surfaces.

\section{II-flat PT Surfaces}

In this section we deal with Riemannian surfaces in Euclidean 3-space $\mathbb{E}^3$
having positive definite second fundamental form $II$. One
may associate to such a surface $M$ geometrical objects measured by means of its second
fundamental form, as {\em second mean curvature} $H_{II}$
and {\em second Gaussian curvature} $K_{II}$, respectively.

\medskip
Further we study the polynomial translation surfaces (PT surfaces in short) with vanishing
second Gaussian curvature in $\mathbb{E}^3$.
By definition (see \cite{kn:MN08}), a PT surface $M$ is given by an immersion
\begin{equation}
\label{mim:eq1}
r:\mathrm{U}\subseteq \mathbb{R}^2\rightarrow \mathbb{R}^3,
\ \ \ \ r(u,v)=(u,v, f(u)+g(v))
\end{equation}
where $f$ and $g$ are polynomial functions.

The first fundamental form $I$ of the surface $M$ is defined by

$$I=Edu^2+2Fdudv+Gdv^2$$\\ where $E=<r_u,r_u>$, $F=<r_u,r_v>$, $G=<r_v,r_v>$ are
the coefficients of $I$, with $<,>$ denoting the scalar product
in $\mathbb{E}^3$ and $r_u=\frac{\partial r(u,v)}{\partial u}$.

The second fundamental form $II$  of $M$ is given by
$$II=edu^2+2fdudv+gdv^2$$\\ where $e=\frac{(r_u,r_v,r_{uu})}{\sqrt{EG-F^2}}$,
$f=\frac{(r_u,r_v,r_{uv})}{\sqrt{EG-F^2}}$ and
$g=\frac{(r_u,r_v,r_{vv})}{\sqrt{EG-F^2}}$.

\medskip
In our case, denoting $f'$ by $\alpha$ and $g'$ by $\beta$, the first fundamental
form $\mathbf{I}$ and the second fundamental form
$\mathbf{II}$ have the following expressions:
$$
{\bf{I}}=\left(1+\alpha(u)^2\right)du^2+2 \alpha(u)\beta(v)du\ dv+\left(1+\beta(v)^2\right)dv^2
$$
$$
{\bf II}=\frac1{\sqrt\Delta}\ \big(\alpha '(u)\ du^2+\beta '(v)\ dv^2\big)
$$
where $\Delta=1+\alpha (u)^2+\beta (v)^2$.

A this point one can immediately compute the Gaussian and the mean curvature using the formulas
$K=\frac{eg-f^2}{EG-F^2}$ and $H=\frac{Eg-2Ff+Ge}{EF-F^2}$ and gets in the case of translation surfaces

\begin{equation}
\label{mim:eq2} \displaystyle
K=\frac{\alpha'(u)\beta'(v)}{(1+\alpha^2+\beta^2)^2} \hspace{5mm}\rm and
\hspace{5mm}
\end{equation}
\begin{equation}
\label{mim:eq3}
\displaystyle H=\frac{(1+\beta^2)\alpha'(u) + (1+\alpha^2)\beta'(v)}{2(1+\alpha^2+\beta^2)^{3/2}}.
\end{equation}

Furthermore, using the formula of the second Gaussian curvature, a similar one to
Brioschi's formula for the Gaussian curvature, obtained replacing the components of the
first fundamental form $E$, $F$, $G$ by those of
the second fundamental form $e$, $f$, $g$

{\small
$$K_{II}=\frac{1}{(|eg|-f^2)^2}\left( \left|
\begin{array}{l}
-\frac{1}{2}e_{vv}+f_{uv}-\frac{1}{2}g_{uu}\quad \frac{1}{2}e_{u} \quad f_{u}-\frac{1}{2}e_{v}\\[2mm]
f_{v} - \frac{1}{2}g_{u} \hspace{22mm} e \hspace{12mm} f \\[2mm]
\frac{1}{2}g_{v}\hspace{28mm} f \hspace{14mm} g
\end{array}
 \right|
 - \left\vert
\begin{array}{l}
0 \hspace{6mm} \frac{1}{2}e_{v} \hspace{5mm} \frac{1}{2}g_{u}\\[2mm]
\frac{1}{2}e_{v}\hspace{5mm} e \hspace{7mm} f\\[2mm]
\frac{1}{2}g_{u}\hspace{5mm} f \hspace{7mm} g
\end{array}
\right\vert \right)$$ }

one gets for the translation surfaces with previous notations

$$K_{II}=\frac{num}{4\alpha'\beta'\Delta^{3/2}},$$

where
\begin{equation}
\label{mim:eq4}
\begin{array}{c}
num=-2\alpha(u)^2\alpha'(u)^2\beta'(v)-2\alpha'(u)\beta(v)^2\beta'(v)^2+\\
2\alpha(u)^2\alpha'(u)\beta'(v)^2 + 2\alpha'(u)^2\beta(v)^2\beta'(v)+\\
2\alpha'(u)\beta'(v)^2+2\alpha'(u)^2\beta'(v) +\\
\alpha'(u)\beta(v)\beta''(v)+\alpha(u)\alpha''(u)\beta'(v)+\\
\alpha(u)^2\alpha'(u)\beta(v)\beta''(v)+\alpha(u)\alpha''(u)\beta(v)^2\beta'(v)+\\
 \alpha'(u)\beta(v)^3\beta''(v)+\alpha(u)^3\alpha''(u)\beta'(v).
\end{array}
\end{equation}

\medskip
It turns out that we have to find those polynomials
$\alpha$ and $\beta$ of degree $m$, respectively $n$, so that $num =
0$ in \eqref{mim:eq4}. At this point let us consider

$$
\alpha=a_mu^m+a_{m-1}u^{m-1}+\ldots +a_{1}u+a_{0}\quad {\rm and}\quad
\beta=b_nv^n+b_{n-1}v^{n-1}+\ldots+b_{1}v+b_{0}
$$
where $a_m$ and $b_n$ are different from $0$.
Replacing $\alpha$ and $\beta$ in \eqref{mim:eq4} we obtain
a polynomial expression in $u$ and $v$ vanishing identically.
This means that all the coefficients are 0.

\bigskip

Let us distinguish several cases:

\medskip

{\bf Case 1:} $m,n\geq 2$, i.e. $\alpha''\neq 0$ and $\beta''\neq 0$.

\medskip

{\bf a.} Suppose $m>n (\geq2)$

The dominant term corresponds to $u^{4m-2}v^{n-1}$ and it comes from\\
$-2\alpha^2\alpha '\beta '+\alpha^3\alpha ''\beta ',$ having the
coefficient $a_{m}^4b_{n}mn(-m-1)$. This expression cannot vanish since $a_m,b_n\neq 0$ and $m>n\geq2$.

{\bf b.} Suppose $n>m (\geq2)$\\
This case can be treated in similar way.

{\bf c.} Suppose $m=n (\geq2)$\\
Analogously, this case cannot occur.

\smallskip
{\bf Case 2:} $m>n=1$, i.e. $\beta=av+b\ {\rm with}\  a,b\in \mathbb{R}\ {\rm and}\ a\neq 0$.

We rewrite the condition $num=0$ in \eqref{mim:eq4} in the following way
\begin{equation}
\label{mim:eq5}
\begin{array}{c}
-2a\alpha(u)^2\alpha'(u)^{2}-2a^2\alpha'(u)\beta(v)^2+2a^2\alpha(u)^2\alpha'(u)\\
+2a\alpha'(u)^2\beta(v)^2+2a^2\alpha'(u)+2a\alpha'(u)^2+\\
a\alpha(u)\alpha''(u)+a\alpha(u)\alpha''(u)\beta(v)^2+a\alpha(u)^3\alpha''(u)=0.
\end{array}
\end{equation}

Using the same idea like in {\bf Case 1}, we analyze the terms of maximum degree in $u$,
namely $u^{4m-2}$, which comes from the expression
$-2a\alpha(u)^2\alpha'(u)^{2}+\alpha(u)^3\alpha''(u)\beta'(v)$ in \eqref{mim:eq5}
having the coefficient $a_m^4am(-m-1)$. This can not vanish
since $a_m^4, a \neq 0$ and $m>n=1$.

\smallskip
The subcase $m=n=1$, treated in similar way, still can not occur.

\smallskip
{\bf Case 3:} $m\geq n=0$ (or $n\geq m=0$)
If $\alpha$ (or $\beta$) is a polynomial of degree $0$, namely is constant,
the condition $num=0$ is automatically satisfied. The parametrization
of the PT surface with vanishing second Gaussian curvature can be written in one of the following forms

\begin{eqnarray}
\label{mim:eq6}
r(u,v)=\big( u,v, au+g(v) \big)\\
\label{mim:eq7}
r(u,v)=\big( u,v, f(u)+a v \big)
\end{eqnarray}
where $f$ and $g$ are arbitrary polynomials, $a\in\R$ (it can also vanishes).
These two surfaces  (given by \eqref{mim:eq6} and \eqref{mim:eq7}) are both
cylinders. Recall that if the second fundamental form for a surface is degenerated it is not possible to
define the second Gaussian curvature.
\smallskip

At this point we formulate the following non-existence result

\begin{theorem}
There are not II-flat polynomial translation surfaces {\rm (i.e. with vanishing second Gaussian curvature)}.
\end{theorem}

\medskip

Until now we studied PT surfaces. Inspired from the example given by Blair
in \cite{kn:BK92} and mentioned in {\em Introduction}, we
analyze other types of translation surfaces, involving {\em power functions}, i.e.
\begin{center}
$\alpha = au^p$ and $\beta=bv^q$ with $a,b\in \mathbb{R},\ a,b\neq
0$ and $p,q\in \mathbb{Q}$.
\end{center}

Remark that because we deal with rational numbers $p$ and $q$, we have to restrict
the coordinate functions to be positive.

The condition for vanishing second Gaussian curvature becomes:

\begin{equation}
\label{mim:eq8}
\begin{array}{c}
a(3p-1)u^pv+a^2b(3q-1)u^{2p+1}v^q + a^3(-p-1)u^{3p}v+\\
b(3q-1)uv^q+ab^2(3p-1)u^pv^{2q+1}+b^3(-q-1)uv^{3q}=0
\end{array}
\end{equation}

\medskip
Again, using the same technique proposed before in this article, we analyze
expression \eqref{mim:eq8} and the only possible case is for
degrees $p=q=\frac{1}{3}$, which implies the additional condition for coefficients,
$a=-b$. We get the surface $M$ given by the parametrization
\begin{equation}
\label{mim:eq9}
r(u,v)=\left( u,v,c(u^{\frac 43}-v^{\frac 43})\right), c\in\R.
\end{equation}
We remark that, up to the multiplication factor $c$, the example given by Blair is the
only one translation surface of this type with
vanishing second Gaussian curvature.

\begin{proposition}
The only translation surfaces given by an immersion \eqref{mim:eq1} where $f$ and $g$
are this time power functions with vanishing second
Gaussian curvature can be parametrized by \eqref{mim:eq9}.
\end{proposition}

\section{$\{K_{II}, H\}$ -- Generalized Weingarten translation surfaces}

If $A,B$ are two different type curvatures of a (nondevelopable) surface,
and if there is a non-trivial functional relation between $A$ and $B$,
then the surface is called an $\{A,B\}$ -- generalized Weingarten surface.
See for details \cite{kn:Sod05}.

\medskip

Among all these surfaces we are interested to study those involving
the second Gaussian curvature $K_{II}$ and the mean curvature $H$.

\subsection{${\mathbf{K_{II}=H}}$}\label{ssec31}

As we have already mentioned in the {\em Introduction}, one interesting property
intensively studied in last years for different types of surfaces, is $K_{II}=H$.
In this subsection we investigate the translation surfaces having this property.

\medskip

Using the previous results, from \eqref{mim:eq3} and \eqref{mim:eq4} the relation $$K_{II}=H$$ is equivalent with

\begin{equation}
\label{mim:eq10}
num=2\alpha'\beta'(\alpha'+\beta^2\alpha'+\beta'+\alpha^2\beta').
\end{equation}

We remark that some terms cancel in \eqref{mim:eq10} and then, dividing by $\alpha'\beta'\neq 0$ we get

\begin{equation}
\label{mim:eq11}
-2\alpha^2\alpha'-2\beta^2\beta'+\frac{\beta\beta''}{\beta'}+
\frac{\alpha\alpha''}{\alpha'}+\alpha^2\frac{\beta\beta''}{\beta'}+
\frac{\alpha\alpha''}{\alpha'}\beta^2+\frac{\beta^3\beta''}{\beta'}+\frac{\alpha^3\alpha''}{\alpha'}=0
\end{equation}

Now, after successive derivations with respect to $u$ and $v$ in \eqref{mim:eq11} one gets

$$\frac {\partial}{\partial u}(\alpha^2)\frac {\partial}{\partial v}\left(\frac{\beta\beta''}{\beta'}\right) +
\frac {\partial}{\partial u}
\left(\frac{\alpha\alpha''}{\alpha'}\right)\frac {\partial}{\partial v}(\beta^2)=0
$$
which is equivalent to
\begin{equation}
\label{mim:eq12}
-\frac{1}{2\alpha\alpha'} \frac {\partial}{\partial u} \left(\frac{\alpha\alpha''}{\alpha'}\right)=
\frac{1}{2\beta\beta'}\frac {\partial}{\partial v}\left(\frac{\beta\beta''}{\beta'}\right).
\end{equation}
Observing now that the two members in \eqref{mim:eq12} are functions of $u$, respectively
of $v$, thus the equality holds only in the case in
which both members are real constants, namely:
$$
\frac{1}{2\alpha\alpha'} \frac {\partial}{\partial u}
\left(\frac{\alpha\alpha''}{\alpha'}\right)=-c \hspace{10mm} {\rm and} \hspace{10mm}
\frac{1}{2\beta\beta'}\frac {\partial}{\partial v}\left(\frac{\beta\beta''}{\beta'}\right)=c, \ c\in \mathbb{R}.
$$
A first integration yields
$$
\frac{\alpha\alpha''}{\alpha'}=-c\alpha^2+d_1 \hspace{10mm}{\rm and}{\hspace{10mm}}
\frac{\beta\beta''}{\beta'}=c\beta^2+d_2,
$$
where $d_1,d_2\in \mathbb{R}$ are integration constants.
A second integration leads to two ODE's fulfilled by $\alpha$ and $\beta$, namely

\begin{equation}
\label{mim:eq13}
2\alpha'=-c\alpha^2+2d_1\ln |\alpha |+ 2m_1 \hspace{3mm}{\rm and}{\hspace{3mm}}
2\beta'=c\beta^2+2d_2\ln | \beta | +2m_2, \ m_1,m_2\in \mathbb{R}.
\end{equation}

In order to solve these equations we distinguish the following cases concerning
the real constants involved in the previous expressions.

\vspace{5mm}
{\bf Case 1.} First, let's suppose that $c\neq 0$. Replacing the expressions from
\eqref{mim:eq13} in \eqref{mim:eq11} we obtain that $\alpha$ and $\beta$ must satisfy
\begin{equation}
\label{mim:eq14}
p_1\alpha^2 - 2d_1\alpha^2\ln |\alpha| + p_2\beta^2 - 2d_2\beta^2\ln |\beta| + d_1 + d_2 = 0
\end{equation}
where $p_1=-c-2m_1+d_1+d_2$ and $p_2=c-2m_2+d_1+d_2$.

\medskip

We remark that in \eqref{mim:eq14} we have a sum of two functions depending on
$u$ respectively on $v$, hence they are constants having opposite signs.
More precisely, the following relations hold
\begin{equation}
\label{mim:eq15}
p_1\alpha^2 -2d_1\alpha^2\ln |\alpha| + d_1 = -q \hspace{5mm} {\rm and}\hspace{5mm}
p_2\beta^2-2d_2\beta^2\ln |\beta| + d_2 = q.
\end{equation}

An elementary study of the function $f:(0,\infty)\rightarrow \mathbb{R}$,
$f(x)=px^2-2dx^2\ln x + d$, where $p,d\in\mathbb{R}$ with $p^2+d^2\neq 0$ shows us that both
equations in \eqref{mim:eq15} have at most two solutions, which means that
$\alpha$ and $\beta$ should be real constants, but this is in contradiction
with $\alpha'\beta'\neq 0$.

If $p_1=d_1=0$ and $p_2=d_2=0$, then equation \eqref{mim:eq14} becomes identity.
But this occurs when $m_1=-\frac{c}{2}$ and $m_2=\frac{c}{2}$.
So, the equations
\eqref{mim:eq13} have the following solutions

$$\alpha(u)=-\tan \frac {cu}{2} \hspace{10mm} {\rm and}\hspace{10mm}
\beta(v)=\tan \frac{cv}{2},$$
considering the integration constants equal to zero, fact which is equivalent
to a possible translation of the parameters.

\vspace{5mm}
{\bf Case 2.} Let's suppose now that $c=0$. In the same manner, by straightforward
computations we get again that $\alpha$ and $\beta$  should be real constants.

\vspace{5mm}
Recall the fact that, at a certain moment, we divided by $\alpha'\beta'\neq 0$.

\smallskip
We ask now what happens when $\alpha'=0$ or $\beta'=0$? We immediately obtain that
$K_{II}=0$. Having in mind $K_{II}=H$ we get $\alpha'=\beta'=0$. The surface is parametrized by
$$
r(u,v)=(u,v, au+bv+c),\ a,b,c\in\mathbb{R},
$$
which is a (portion of a) plane, having degenerated second fundamental form. So, this case cannot occur.

Now, we can state the following theorem:

\begin{theorem}
The only translation surfaces with non-degenerate second fundamental form having the property $K_{II}=H$
are given, up to a rigid motion of $\mathbb{R}^3$, by

\begin{equation}
\label{mim:eq16}
\displaystyle
r(u,v)=\left(u,\ v,\ \frac{2}{c}\ \ln\left| \frac{\cos \frac{cu}{2}}{\cos\frac{cv}{2}}\right| \ \right), \ c\in \mathbb{R}^{*}.
\end{equation}

More, we remark that \eqref{mim:eq16} is the parametrization of a Scherk type surface, so we have $$K_{II}=H=0.$$
\end{theorem}

\vspace{5mm}

\subsection{${\mathbf{K_{II}=\lambda H, \ \lambda\neq 1,2}}$}
We are wondering what happens in the general case
\begin{equation}
\label{mim:eq17}
K_{II}=\lambda H,\ \lambda\neq 1
\end{equation}
for translation surfaces given by parametrization \eqref{mim:eq1}.
\smallskip

Following the same idea as in the previous case for $\lambda=1$,
using the formulae \eqref{mim:eq3} and \eqref{mim:eq4} the relation \eqref{mim:eq17}
becomes equivalent with

\begin{equation}
\label{mim:eq18}
num = 2\lambda \alpha'\beta'(\alpha'+\beta^2\alpha'+\beta'+\alpha^2\beta').
\end{equation}

Dividing by $\alpha'\beta'$ \eqref{mim:eq18} implies that

\begin{equation}
\label{mim:eq19}
\displaystyle
\begin{array}{c}
\displaystyle
-2\alpha^2\alpha'-2\beta^2\beta'+ \frac{\alpha^3\alpha''}{\alpha'}+\frac{\beta^3\beta''}{\beta'}+\\[3mm]
\displaystyle
\left((2-2\lambda)\alpha'+\frac{\alpha\alpha''}{\alpha'}\right)(\beta^2+1)+
\left((2-2\lambda)\beta'+\frac{\beta\beta''}{\beta'}\right)(\alpha^2+1)=0.
\end{array}
\end{equation}

After successive derivations with respect to $u$ and $v$ in \eqref{mim:eq19} one gets

$$
\displaystyle
-\frac{1}{2\alpha\alpha'} \frac {\partial}{\partial u} \left(2(1-\lambda)\alpha'+ \frac{\alpha\alpha''}{\alpha'}\right)=
\frac{1}{2\beta\beta'}\frac {\partial}{\partial v}\left(2(1-\lambda)\beta'+\frac{\beta\beta''}{\beta'}\right).
$$

Using the same technique as in subsection \ref{ssec31}, notice that the previous equality
can occur only when both sides are equal with the same real constant denoted by $\mu$, that means we have
$$
\frac{1}{2\alpha\alpha'} \frac {\partial}{\partial u}
\left(2(1-\lambda)\alpha'+\frac{\alpha\alpha''}{\alpha'}\right)=-\mu
\hspace{2mm} {\rm and} \hspace{2mm}
\frac{1}{2\beta\beta'}\frac {\partial}{\partial v}\left(2(1-\lambda)\beta'+
\frac{\beta\beta''}{\beta'}\right)=\mu.
$$

After a first integration one obtains
\begin{equation}
\label{mim:eq20}
2(1-\lambda)\alpha'+\frac{\alpha\alpha''}{\alpha'}=-\mu\alpha^2+\nu_1 \hspace{5mm}{\rm and}{\hspace{5mm}}
2(1-\lambda)\beta'+\frac{\beta\beta''}{\beta'}=\mu\beta^2+\nu_2,
\end{equation}
where $\nu_1$ and $\nu_2$ are integration constants.
\medskip

Then, a second integration yields two ODE's fulfilled by $\alpha$ and $\beta$,
analogously with equations \eqref{mim:eq13}, namely
\begin{equation}
\label{mim:eq21}
\begin{array}{c}
\displaystyle
\alpha'=\tau_1\alpha^{2\lambda-2}-\frac{\mu}{4-2\lambda}\alpha^2+\frac{\nu_1}{2-2\lambda}\\[3mm]
\displaystyle
\beta'=\tau_2\beta^{2\lambda-2}+\frac{\mu}{4-2\lambda}\beta^2+\frac{\nu_2}{2-2\lambda}
\end{array}
\end{equation}
where $\tau_1$ and $\tau_2$ are integration constants. Clearly, $\lambda\neq 2$.

\smallskip
In the most general case, considering that all the integration constants are
different form zero and replacing \eqref{mim:eq20} and \eqref{mim:eq21} in
\eqref{mim:eq19} we get the following expressions satisfied by $\alpha$ and $\beta$:

\begin{equation}
\label{mim:eq22}
\begin{array}{c}
\displaystyle
(2\lambda-4)\tau_1\alpha^{2\lambda}+\left(-\mu-\frac{\nu_1}{1-\lambda}+\nu_2\right)\alpha^2+\nu_1+\qquad\\
\displaystyle
\qquad +(2\lambda-4)\tau_2\beta^{2\lambda}+\left(\mu+\nu_1-\frac{\nu_2}{1-\lambda}\right)\beta^2+\nu_2=0.
\end{array}
\end{equation}

It is easy to remark that \eqref{mim:eq22} is a sum of two functions, one depending on $u$,
the second depending on $v$. Consequently, the equality holds if these
functions are respectively plus and minus a constant. More, at this point, the algebraic
equations we obtained, having the coefficients different from zero, have a finite number of solutions.
This means that $\alpha$ and $\beta$ are constants, case in which the second
fundamental form is degenerate. So the situation cannot occur.

\smallskip
Further, we analyze the case when all the coefficients in the two algebraic equations are identically zero,
namely
\begin{equation}
\label{mim:eq23}
 \left\{
\begin{array}{ll}
    (2\lambda-4)\tau_1=0 \quad &
\displaystyle
    -\mu-\frac{\nu_1}{1-\lambda}+\nu_2=0 \\[2mm]
    (2\lambda-4)\tau_2=0 &
\displaystyle
    \mu+\nu_1-\frac{\nu_2}{1-\lambda}=0\\[2mm]
    \nu_1+\nu_2=0.
\end{array}
\right.
\end{equation}

Recall that $\lambda\in \mathbb{R}\setminus \{1,2\}$, so \eqref{mim:eq23} yields
\begin{equation}
\label{mim:eq24}
\begin{array}{c}
    \tau_1=\tau_2=0\qquad
    \nu_1=-\nu_2 \qquad
    \displaystyle
    \mu=\frac{\nu_1(\lambda-2)}{1-\lambda}\ .
\end{array}
\end{equation}

Replacing these expressions in \eqref{mim:eq21} we get

$$
\alpha'=\frac{\nu_1}{2-2\lambda}(\alpha^2+1) \hspace{5mm} {\rm and\ } \hspace{5mm}
\beta'=-\frac{\nu_1}{2-2\lambda}(\beta^2+1).
$$

\vspace{5mm}
We can state the following theorem.

\begin{theorem}
The only $\{K_{II},H\}$--generalized Weingarten translation surfaces with non-degen\-erate second fundamental form
satisfying $K_{II}=\lambda H$ with $\lambda\in\mathbb{R}\setminus\{1,2\}$, are given,
up to a rigid motion of $\mathbb{R}^3$, by the parametrization
\begin{equation}
\label{mim:eq25}
r(u,v)=\left(u,v, \frac{1}{p}\log\left|\frac{\cos(pv+r)}{\cos(pu+q)}\right|\right),\
{\rm where\ } p\neq 0  {\rm \ and\ } r,q\in\mathbb{R}
\end{equation}
which represents a Scherk type surface. Moreover $K_{II}=H=0$.
\end{theorem}
\proof
After a straightforward computation we conclude with $p=\frac{\nu_1}{2-2\lambda}\neq 0$,
since $\nu_1$ cannot be zero. Otherwise, if $\nu_1=0$, then $\nu_2=\mu=0$ and the second
fundamental form would be degenerated, which is false.

\endproof

\vspace{5mm}

\subsection{${\mathbf{K_{II}=2H}}$}

In this subsection we deal with $\lambda=2$.

Making similar computations as in the general case, one gets that \eqref{mim:eq19} is equivalent with
\begin{equation}
\label{mim:eq26}
\displaystyle
\begin{array}{c}
\displaystyle
\left(\alpha^2+\beta^2+1\right)\left(-2\alpha'+\frac{\alpha\alpha''}{\alpha'}-2\beta'+\frac{\beta\beta''}{\beta'}\right)=0.
\end{array}
\end{equation}
Hence, using the same strategy as so far, there exists $\nu\in{\mathbb{R}}$ such that $-2\alpha'+\frac{\alpha\alpha''}{\alpha'}=\nu$ and
$-2\beta'+\frac{\beta\beta''}{\beta'}=-\nu$, which yield, after a first integration,
\begin{equation}
\label{mim:eq27}
\alpha'=\tau_1\alpha^2-\frac{\nu}{2}\ \ {\rm and} \ \ \beta'=\tau_2\beta^2+\frac{\nu}{2}
\end{equation}
where $\tau_1$, $\tau_2$ are real constants.

\vspace{5mm}
In order to solve these ODE's we analyze the following cases:

\vspace{5mm}
{\bf Case 1.} $\tau_1, \tau_2\neq 0.$

Rewrite previous relations in the form
$\alpha'=\tau_1\left(\alpha^2-\frac{\nu}{2\tau_1}\right)$  and
$\beta'=\tau_2\left(\beta^2+\frac{\nu}{2\tau_2}\right)$.

\medskip

At this point we distinguish some different cases involving the signs of
$\displaystyle \frac{\nu}{\tau_1}$ and $\displaystyle \frac{\nu}{\tau_2}$.

\medskip

{\bf Subcase 1.1} \ $\displaystyle \frac{\nu}{\tau_1}>0$, $\ \displaystyle \frac{\nu}{\tau_2}>0$.

\vspace{5mm}
Denoting $\displaystyle \frac{\nu}{2\tau_1}=\rho^2$,
$\ \displaystyle \frac{\nu}{2\tau_2}=\eta^2$, ($\rho$, $\eta>0$), equations \eqref{mim:eq27} are equivalent with

$$
\frac{\alpha'}{\alpha^2-\rho^2}=\tau_1\ \ \ {\rm and}\ \ \ \frac{\beta'}{\beta^2+\eta^2}=\tau_2.
$$

Integrating, we get two solutions for $\alpha$, namely
\begin{center}
$\alpha=-\rho\coth(\rho\tau_1 u)$, respectively
$\alpha=-\rho\tanh(\rho\tau_1 u)$,  and
$\beta=\eta\tan(\eta\tau_2 v)$,
\end{center}
by considering the integration constants equal to zero.
Finally, we obtain the functions
$f$ and $g$ from the parametrization \eqref{mim:eq1}
$$
f(u)=-\frac{1}{\tau_1}\log(\sinh(\rho\tau_1 u)), {\rm respectively\ } \ f(u)=-\frac{1}{\tau_1}\log(\cosh(\rho\tau_1 u))
$$
$$
{\rm and}\ \ \ g(v)=-\frac{1}{\tau_2}\log(\cos(\eta\tau_2 v)).
$$

Denoting $p:=\rho\tau_1$ and $q:=\eta\tau_2$, we conclude this subcase with two type of translation surfaces, namely
\begin{equation}
\label{mim:eq28}
\displaystyle
r(u,v)=\left(u,v, -\frac{\nu}{2}\log\left(\sinh(pu)^{\frac{1}{p^2}}\cos(qv)^{\frac{1}{q^2}}\right) \right),
\end{equation}

\begin{equation}
\label{mim:eq29}
\displaystyle
r(u,v)=\left(u,v, -\frac{\nu}{2}\log\left(\cosh(pu)^{\frac{1}{p^2}}\cos(qv)^{\frac{1}{q^2}}\right) \right).
\end{equation}

\vspace{5mm}

Note that for all the rest of the situations, the coefficients $p$ and $q$ have the same significance.
The other situations can be treated in similar way, and we will relieve only the main results.

\vspace{5mm}

{\bf Subcase 1.2} \ $\displaystyle \frac{\nu}{\tau_1}<0$, $\ \displaystyle \frac{\nu}{\tau_2}>0$.

\vspace{5mm}
Denoting $\displaystyle -\frac{\nu}{2\tau_1}=\rho^2$, $\ \displaystyle \frac{\nu}{2\tau_2}=\eta^2$,
($\rho$, $\eta>0$), equations \eqref{mim:eq27} are equivalent with
$$
\frac{\alpha'}{\alpha^2+\rho^2}=\tau_1\ \ \ {\rm and}\ \ \ \frac{\beta'}{\beta^2+\eta^2}=\tau_2.
$$

Integrating two times we obtain the parametrization \eqref{mim:eq1}
\begin{equation}
\label{mim:eq30}
\displaystyle
r(u,v)=\left(u,v, \frac{\nu}{2}\ \log\frac{\cos(pu)^{\frac{1}{p^2}}}{\cos(qv)^{\frac{1}{q^2}}} \right).
\end{equation}

\vspace{8mm}

{\bf Subcase 1.3} \ $\displaystyle \frac{\nu}{\tau_1}>0$, $\ \displaystyle \frac{\nu}{\tau_2}<0$.

\vspace{5mm}
Denoting $\displaystyle \frac{\nu}{2\tau_1}=\rho^2$, $\ \displaystyle -\frac{\nu}{2\tau_2}=\eta^2$,
($\rho$, $\eta>0$), equations \eqref{mim:eq27} are equivalent with
$$
\frac{\alpha'}{\alpha^2-\rho^2}=\tau_1\ \ \ {\rm and}\ \ \ \frac{\beta'}{\beta^2-\eta^2}=\tau_2.
$$

Performing same steps, the parametrization of the surface is given by

\begin{equation}
\label{mim:eq31}
\displaystyle
r(u,v)=\left(u,v, -\frac{\nu}{2}\ \log\frac{\sinh(pu)^{\frac{1}{p^2}}}{\sinh(qv)^{\frac{1}{q^2}}} \right),
\end{equation}

\begin{equation}
\label{mim:eq32}
\displaystyle
r(u,v)=\left(u,v, -\frac{\nu}{2}\ \log\frac{\cosh(pu)^{\frac{1}{p^2}}}{\cosh(qv)^{\frac{1}{q^2}}} \right),
\end{equation}

\begin{equation}
\label{mim:eq33}
\displaystyle
r(u,v)=\left(u,v, -\frac{\nu}{2}\ \log\frac{\cosh(pu)^{\frac{1}{p^2}}}{\sinh(qv)^{\frac{1}{q^2}}} \right),
\end{equation}

\begin{equation}
\label{mim:eq34}
\displaystyle
r(u,v)=\left(u,v, -\frac{\nu}{2}\ \log\frac{\sinh(pu)^{\frac{1}{p^2}}}{\cosh(qv)^{\frac{1}{q^2}}} \right).
\end{equation}

\medskip

{\bf Subcase 1.4} \ $\displaystyle \frac{\nu}{\tau_1}<0$, $\ \displaystyle \frac{\nu}{\tau_2}<0$.

\vspace{5mm} Denoting $\displaystyle -\frac{\nu}{2\tau_1}=\rho^2$, $\
\displaystyle -\frac{\nu}{2\tau_2}=\eta^2$, ($\rho$, $\eta>0$),
equations \eqref{mim:eq27} are equivalent with
$$
\frac{\alpha'}{\alpha^2+\rho^2}=\tau_1\ \ \ {\rm and}\ \ \ \frac{\beta'}{\beta^2-\eta^2}=\tau_2.
$$

Solving these ODE's and integrating again, the parametrization of the surface is given by
\begin{equation}
\label{mim:eq35}
\displaystyle
r(u,v)=\left(u,v, \frac{\nu}{2}\ \log\left(\cos(pu)^{\frac{1}{p^2}}\sinh(qv)^{\frac{1}{q^2}}\right) \right),
\end{equation}

\begin{equation}
\label{mim:eq36}
\displaystyle
r(u,v)=\left(u,v, \frac{\nu}{2}\ \log\left(\cos(pu)^{\frac{1}{p^2}}\cosh(qv)^{\frac{1}{q^2}}\right) \right).
\end{equation}

Remark that Subcase 1.4 is similar to Subcase 1.1.

\medskip

Throughout this subsection the domains of definition for the parameters $u$ and $v$ are chosen such that
all expressions are well defined, even that we do not mention explicitly this fact.

\vspace{5mm}

{\bf Case 2.} $\tau_1=\tau_2=0$.

Substituting in \eqref{mim:eq27} one gets $\displaystyle\alpha'=-\frac{\nu}{2}$ and
$\displaystyle\beta'=\frac{\nu}{2}$. Integrating twice, and
re-denoting the constants we obtain the hyperbolic paraboloid
\begin{equation}
\label{mim:eq37}
r(u,v)=(u,v, a(u-u_0)^2-a(v-v_0)^2), \ \ a, u_0, v_0\in\mathbb{R}.
\end{equation}

We can state the main theorem in the case $K_{II}=2 H$.

\begin{theorem}
The only translation surfaces with non-degenerate second fundamental form having the property $K_{II}=2 H$
are given, up to a rigid motion of $\mathbb{R}^3$, by
\begin{itemize}
\item [\it i)] parametrizations \eqref{mim:eq28} -- \eqref{mim:eq36};
\item[\it ii)] parametrization \eqref{mim:eq37} of a hyperbolic paraboloid.
\end{itemize}
\end{theorem}
\proof
As we have already seen item (i) occurs if $\tau_1,\tau_2\neq 0$.
The second item holds if both $\tau_1$ and $\tau_2$ vanish.
\endproof


\medskip

Concerning the case of polynomial translation surfaces we conclude with
\begin{theorem}
The only PT surfaces satisfying the relation $K_{II}=\lambda H$ are given by the hyperbolic
paraboloids parametrized by \eqref{mim:eq37}, case in which $\lambda=2$.
\end{theorem}
\proof
The statement is a consequence of the previous theorems.
\endproof

\section{II-minimal translation surfaces}

Similar to the variational characterization of the mean curvature $H$,
the curvature of the second fundamental form, denoted by $H_{II}$ is
introduced as a measure for the rate of change of the II-area under a normal deformation.
For details see \cite{kn:HVV08}.
In this section we analyze II-minimal translation surfaces with non-degenerate second fundamental form,
namely we study under which conditions the second mean curvature vanishes,
i.e. $H_{II}=0$. Having in mind the usual technique for computing the second mean curvature by using the
normal variation of the area functional one gets
$$
H_{II}=-H-\frac{1}{4}\Delta^{II}\log(K)
$$
where $\Delta^{II}$ is the Laplacian for functions computed with respect to the second
fundamental form as metric. $H_{II}$ can be equivalently expressed as
\begin{equation}
\label{mim:eq38}
H_{II}=-H-\frac{1}{2\sqrt{\det II}}\sum \limits_{i,j} \frac{\partial}{\partial u^i}
\left(\sqrt{\det II}\ h^{ij}\ \frac{\partial}{\partial u^j}(\ln \sqrt{K}) \right).
\end{equation}
Here $II$ denotes the second fundamental form defined in {\em Section $2$},
$(h_{ij})$ is the associated matrix with its inverse $(h^{ij})$,
the indices $i,j$ belong to $\{1,2\}$ and the parameters $u^1$, $u^2$ are $u$,
respectively $v$ from the parametrization \eqref{mim:eq1}.
Moreover, $II$ becomes a metric on the surface if it is non-degenerated.
The inverse matrix $(h^{ij})$ has the following expression
$$
\left(h^{ij}\right)_{i,j}=\left(
\begin{array}{l}
\frac{\sqrt{1+\alpha^2+\beta^2}}{\alpha'} \hspace{10mm} 0\\[2mm]
\hspace{5mm}0 \hspace{10mm} \frac{\sqrt{1+\alpha^2+\beta^2}}{\beta'}
\end{array}
\right)
$$
where $K$ and $H$ denote the usual Gaussian, respectively mean curvatures
of our surface.

After straightforward computations, the sum in \eqref{mim:eq38} becomes:
$$
\begin{array}{rcl}
\displaystyle\sum\limits_{i,j} & = & \frac{1}{4\Delta^2}\sqrt{\frac{\beta'}{\alpha'}}
     \left(\frac{2\alpha'\alpha'''-3\alpha''^2}{\alpha'^2}\ \Delta^2
     +(-4\alpha\alpha''-8\alpha'^2)\Delta+16 \alpha^2\alpha'^2\right) + \\[4mm]
& &  +\frac{1}{4\Delta^2}\sqrt{\frac{\alpha'}{\beta'}}
     \left(\frac{2\beta'\beta'''-3\beta''^2}{\beta'^2}\ \Delta^2
     +(-4\beta\beta''-8\beta'^2)\Delta+16\beta^2\beta'^2\right).
\end{array}
$$
We are interested to find $II$--minimal translation surfaces in the
Euclidean 3-space. Having now all the necessary tools, the condition $H_{II}=0$ for a
translation surface is equivalent to
\begin{equation}
\label{mim:eq39}
\frac{2\alpha'\alpha'''-3\alpha''^2}{2\alpha'^3}+\frac{2\beta'\beta'''-3\beta''^2}{2\beta'^3}-
\frac{2}{\Delta} \left(\frac{\alpha'^2+\alpha\alpha''}{\alpha'}+\frac{\beta'^2+\beta\beta''}{\beta'}\right)+
\frac{6}{\Delta^2}\ (\alpha^2\alpha'+\beta^2\beta')=0.
\end{equation}

The first two terms in \eqref{mim:eq39} are functions only of
$u$ respectively of $v$, hence we derive in the previous equation
successively with respect to $u$ and $v$.
\smallskip

Denoting by
$\displaystyle \phi(u)=\frac{\alpha'^2+\alpha\alpha''}{\alpha'}\ $,
$\displaystyle \psi(v)=\frac{\beta'^2+\beta\beta''}{\beta'}\ $,
$p(u)=\alpha^2\alpha'$ and $q(v)=\beta^2\beta'$, we get
$$
\frac{\partial}{\partial v}\frac{\partial}{\partial u}\left(-\frac{2}{\Delta}\ (\phi+\psi)+
\frac{6}{\Delta^2}\ (p+q)\right)=0.
$$
After straightforward computations and multiplying with $\frac{\Delta^3}{8\alpha\alpha'\beta\beta'}$
it follows
\begin{equation}
\label{mim:eq40}
(F+G)\Delta^2-2(P+Q)\Delta+18(p+q)=0,
\end{equation}
where $\displaystyle F(u)=\frac{\phi'}{2\alpha\alpha'}\ $,
$\displaystyle G(v)=\frac{\psi'}{2\beta\beta'}\ $,
$\displaystyle P(u)=\phi+\frac{3p'}{2\alpha\alpha'}$ and
$\displaystyle Q(v)=\psi+\frac{3q'}{2\beta\beta'}\ $.

Repeating the same operations, namely the two partial derivatives and
the division by $4\alpha\alpha'\beta\beta'$ one gets
\begin{equation}
\label{mim:eq41}
(A+B)\Delta+a+b=0,
\end{equation}
where $\displaystyle A(u)=\frac{F'}{2\alpha\alpha'}\ $,
$\displaystyle B(v)=\frac{G'}{2\beta\beta'}\ $,
$\displaystyle a(u)=F-\frac{P'}{2\alpha\alpha'}$ and
$\displaystyle b(v)=G-\frac{Q'}{2\beta\beta'}\ $.

\smallskip

Finally, using the same technique, we should have
$$
\frac{A'}{2\alpha\alpha'}=c \qquad \frac{B'}{2\beta\beta'}=-c,\quad c\in \mathbb{R}.
$$
Solving the above equations we obtain  $A(u)=c\alpha^2+d_1$ and $B(v)=-c\beta^2+d_2$.
Replacing these expressions in the previous ODEs we find that

$\begin{array}{ll}
F(u)=\frac{c}{2}\alpha^4+d_1\alpha^2+\mu_1 \qquad&\qquad
G(v)=-\frac{c}{2}\beta^4+d_2\beta^2+\mu_2\qquad \\[3mm]
\phi(u)=\frac{c}{6}\alpha^6+\frac{d_1}{2}\alpha^4+\mu_1\alpha^2+\tau_1 \qquad &\qquad
\psi(v)=-\frac{c}{6}\beta^6+\frac{d_2}{2}\beta^4+\mu_2\beta^2+\tau_2\qquad
\end{array}$
\begin{equation}
\label{mim:eq42}
\begin{array}{ll}
\alpha'(u)=\frac{c}{42}\alpha^6+\frac{d_1}{10}\alpha^4+\frac{\mu_1}{3}\alpha^2+
        \tau_1+\frac{\rho_1}{\alpha} \quad &
\beta'(v)=-\frac{c}{42}\beta^6+\frac{d_2}{10}\beta^4+
\frac{\mu_2}{3}\beta^2+\tau_2+\frac{\rho_2}{\beta}\quad
\end{array}
\end{equation}
$\begin{array}{ll}
p(u)=\frac{c}{42}\alpha^8+\frac{d_1}{10}\alpha^6+\frac{\mu_1}{3}\alpha^4+
        \tau_1\alpha^2+\rho_1\alpha &
q(v)=-\frac{c}{42}\beta^8+\frac{d_2}{10}\beta^6+\frac{\mu_2}{3}\beta^4+\tau_2\beta^2+\rho_2\beta\\[3mm]
P(u)=\frac{19 c}{42}\alpha^6+\frac{7d_1}{5}\alpha^4+3\mu_1\alpha^2+
        4\tau_1+\frac{3\rho_1}{2\alpha} &
Q(v)=-\frac{19 c}{42}\beta^6+\frac{7d_2}{5}\beta^4+3\mu_2\beta^2+4\tau_2+\frac{3\rho_2}{2\beta}\\[3mm]
a(u)=-\frac{6c}{7}\alpha^4-\frac{9d_1}{5}\alpha^2-2\mu_1+\frac{3\rho_1}{4\alpha^3}&
b(v)=\frac{6c}{7}\beta^4-\frac{9d_2}{5}\beta^2-2\mu_2+\frac{3\rho_2}{4\beta^3}
\end{array}$

where $d_1,d_2,\mu_1,\mu_2,\tau_1,\tau_2,\rho_1,\rho_2 \in \mathbb{R}$.
In order to determine all these integration constants, we substitute
the corresponding expressions in \eqref{mim:eq41}, obtaining a sum of polynomials in $\alpha$ and $\beta$
equals to 0. This means that there exists $\xi\in \mathbb{R}$ such that
$$
\begin{array}{ll}
    \displaystyle \ \ \frac{c}{7}\ \alpha^4+\left(c-\frac{4}{5}\ d_1+d_2\right)\alpha^2+
    \frac{3\rho_1}{4\alpha^3}+d_1-2\mu_1-\xi=0 \\[4mm]
    \displaystyle-\frac{c}{7}\ \beta^4+\left(-c+d_1-\frac{4}{5}\ d_2\right)\beta^2+
\frac{3\rho_2}{4\beta^3}+d_2-2\mu_2+\xi=0 .
\end{array}
$$
At this point, by the same argument as in previous section, all the coefficients in the
above (algebraic) expressions must be zero and consequently we get $c=0$, $d_1=d_2=0$,
$\rho_1=\rho_2=0$, $\mu_1=-\frac{\xi}{2}$ and $\mu_2=\frac{\xi}{2}\ $.
Thus \eqref{mim:eq42} can be expressed in a simpler form
\begin{center}
\begin{equation}
\label{mim:eq43}
\begin{array}{ll}
F(u)=-\frac{\xi}{2} &
G(v)=\frac{\xi}{2}\\[3mm]
\phi(u)=-\frac{\xi}{2}\alpha^2+\tau_1 &
\psi(v)=\frac{\xi}{2}\beta^2+\tau_2\\[3mm]
\alpha'(u)=-\frac{\xi}{6}\alpha^2+\tau_1 &
\beta'(v)=\frac{\xi}{6}\beta^2+\tau_2\\[3mm]
p(u)=-\frac{\xi}{6}\alpha^4+\tau_1\alpha^2 &
q(v)=\frac{\xi}{6}\beta^4+\tau_2\beta^2\\[3mm]
P(u)=-\frac{3\xi}{2}\alpha^2+4\tau_1 &
Q(v)=\frac{3\xi}{2}\beta^2+4\tau_2\\[3mm]
a(u)=\xi &
b(v)=-\xi.
\end{array}
\end{equation}
\end{center}
Let us take a look in \eqref{mim:eq40}. By the same reasoning as above
we deduce
$$
 \begin{array}{ll}
   \ (3\xi+10\tau_1-8\tau_2)\alpha^2-8\tau_1=\eta\\[2mm]
    (-3\xi-8\tau_1+10\tau_2)\beta^2-8\tau_2=-\eta
\end{array}
$$
for an arbitrary $\eta\in \mathbb{R}.$ Moreover,
the integration constants should be $\tau_1=-\frac{\eta}{8}\ $, $\tau_2=\frac{\eta}{8}\ $,
$\xi=\frac{3\eta}{4}$ and we can conclude that $\alpha'=-\frac{\eta}{8}(\alpha^2+1)$ and
$\beta'=\frac{\eta}{8}(\beta^2+1)$.
Finally $\alpha$ and $\beta$ must satisfy also the condition \eqref{mim:eq39}.
After straightforward computations it follows that $\eta=\xi=0$.

The conclusion is $\alpha'=\beta'=0$, which cannot occur since if this happened
the second fundamental form would vanish identically.
Hence the second mean curvature is not defined and we end this section with
the following non-existence theorem
\begin{theorem}
There are not II-minimal translation surfaces in Euclidean 3-space.
\end{theorem}

\begin{appendix}
\begin{center}{\bf Appendix}\end{center}

We conclude this paper with some figures realized with {\em Mathematica $5.0$}
representing surfaces obtained in the classification {\bf Theorem 5}:

\begin{figure}[htb]
\begin{center}
\epsfxsize=37mm \centerline{\leavevmode \epsffile{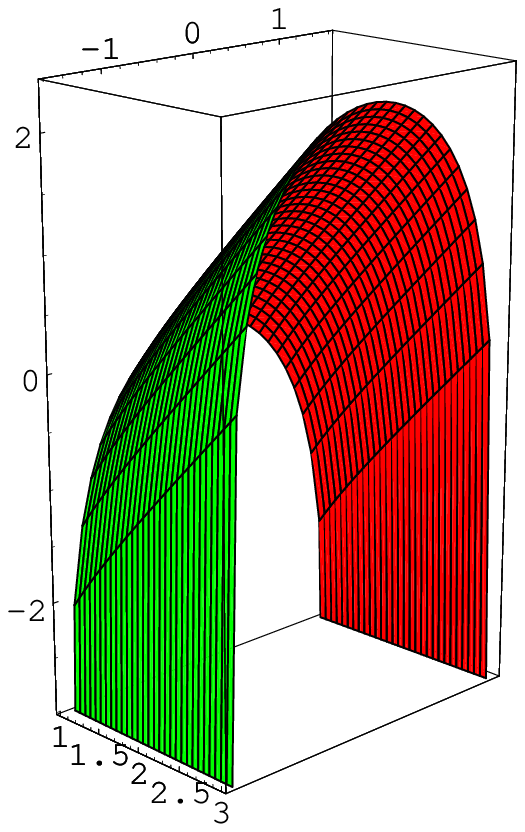} \hspace{20mm}
\epsfxsize=47mm \epsffile{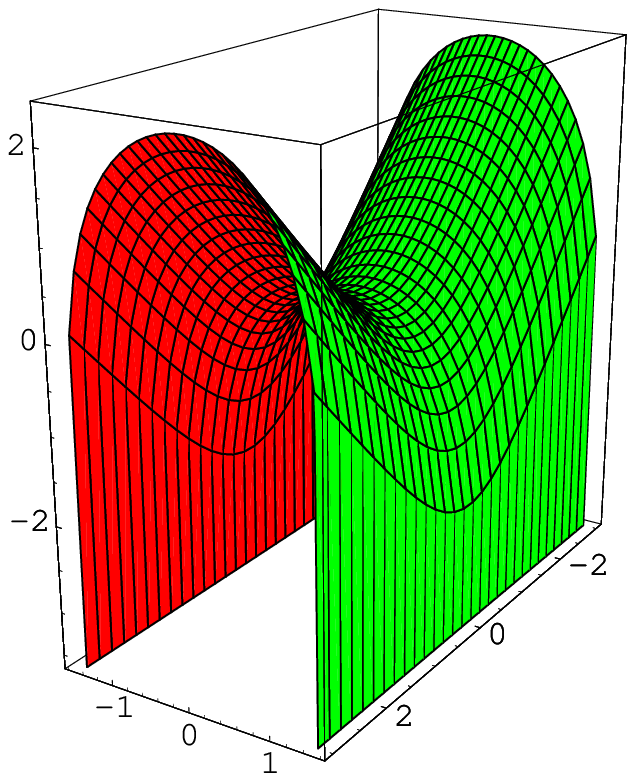} }
\end{center}
\*{$r(u,v)=(u, v,\ln(\sinh u\cos v))$\hspace{20mm} $r(u,v)=(u, v,\ln(\cosh u\cos v))$}
\end{figure}

\pagebreak

\begin{figure}[htb]
\begin{center}
\epsfxsize=65mm \centerline{\leavevmode \epsffile{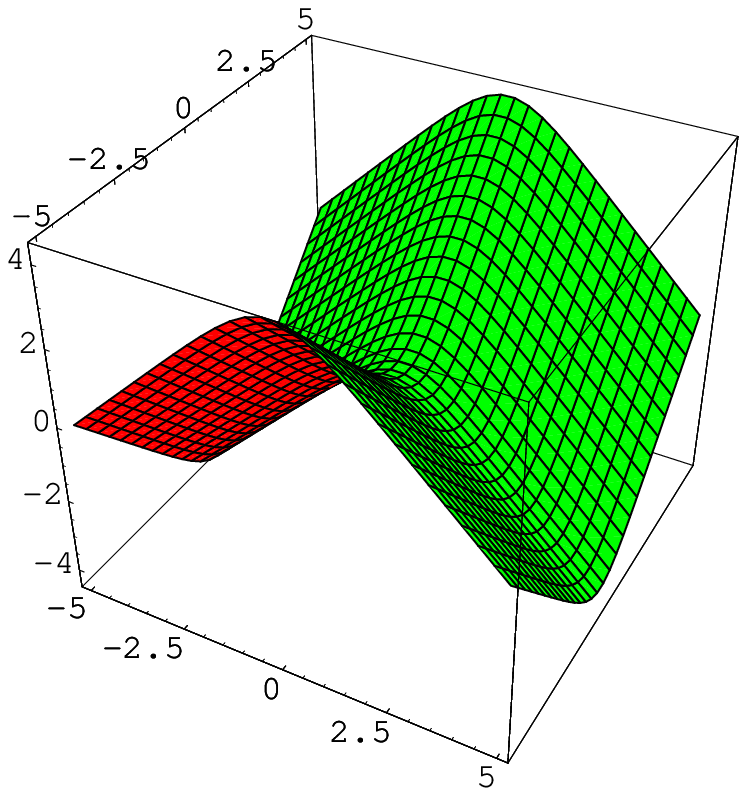} \hspace{20mm}
\epsfxsize=40mm \epsffile{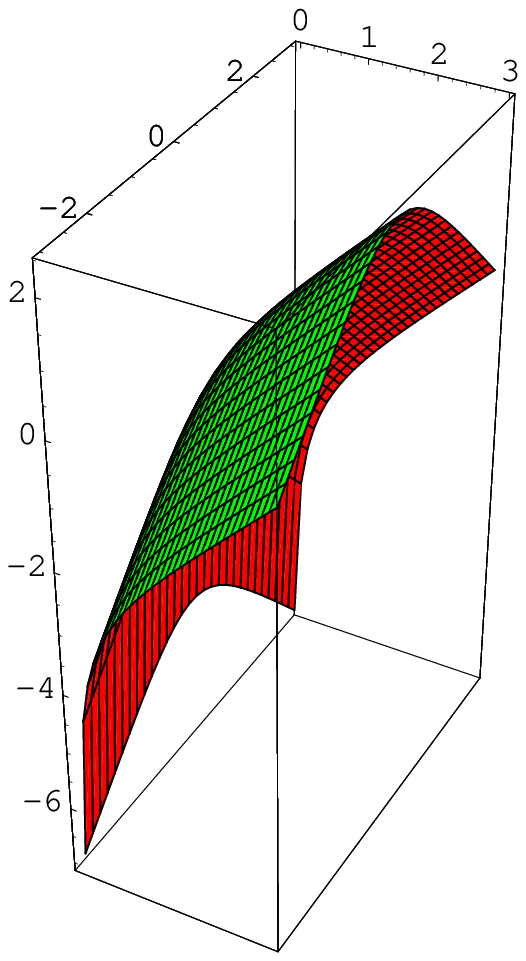} }
\end{center}
\*{$\displaystyle r(u,v)=\left(u,v,\ln\frac{\cosh u}{\cosh v}\right)$\hspace{20mm}
$\displaystyle r(u,v)=\left(u,v,\ln\frac{\sinh u}{\cosh v}\right)$}
\end{figure}

\begin{figure}[htb]
\begin{center}
\epsfxsize=35mm \centerline{\leavevmode \epsffile{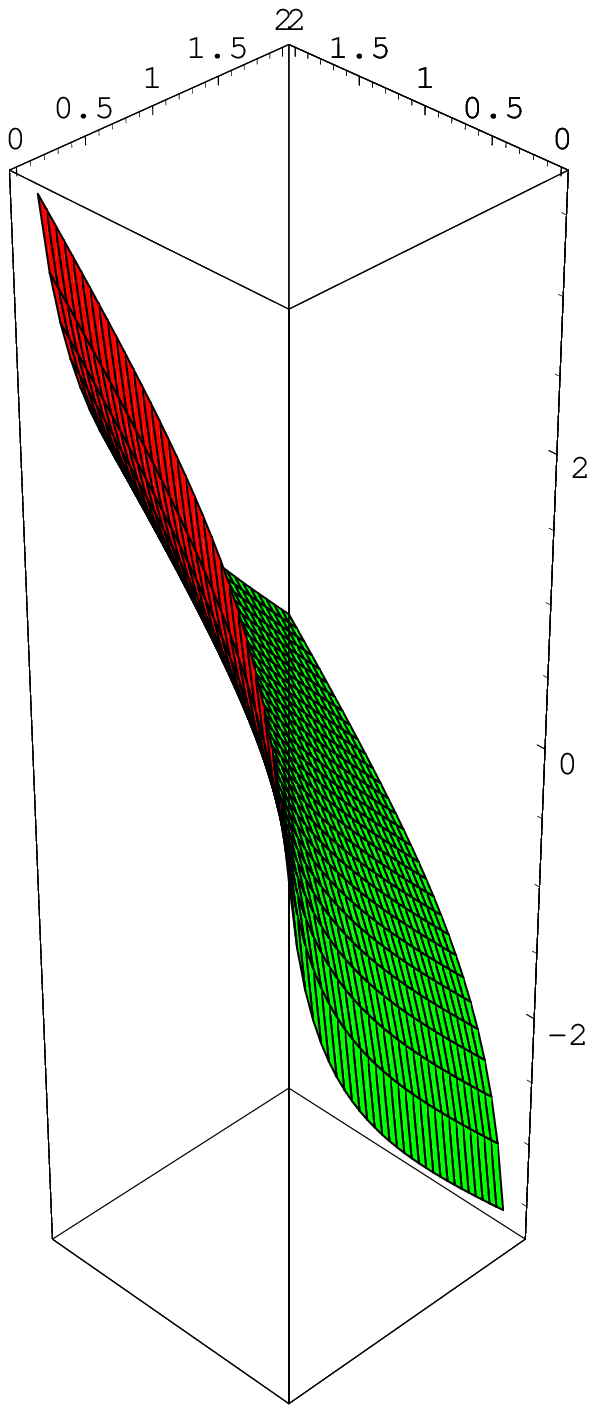} \hspace{30mm}
\epsfxsize=40mm \epsffile{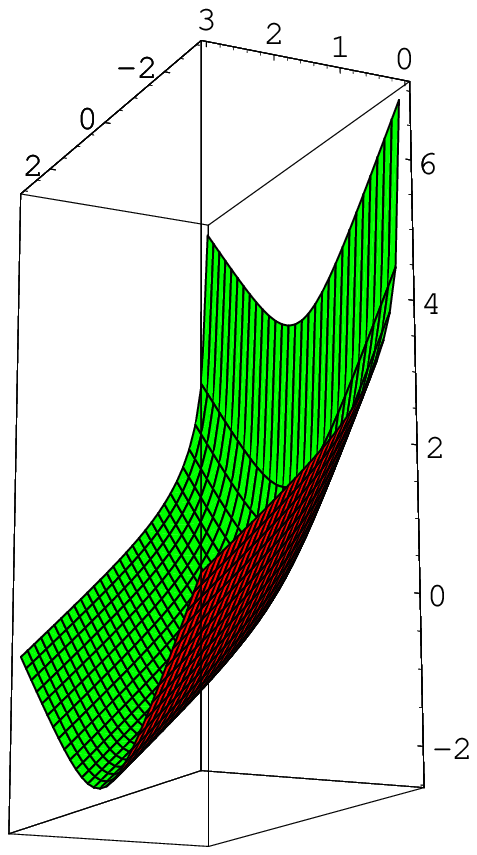} }
\end{center}
\*{$\displaystyle r(u,v)=\left(u,v,\ln\frac{\sinh u}{\sinh v}\right)$\hspace{30mm}
$\displaystyle r(u,v)=\left(u,v,\ln\frac{\cosh u}{\sinh v}\right)$}
\end{figure}

\end{appendix}

{\small {\bf Acknowledgement.} The authors would like to thank the Research Group
of the Department of Geometry, Katholieke Universiteit Leuven, Belgium, in particular
Prof. Franki Dillen for useful comments and critical remarks during the preparation of this paper.}
\pagebreak

\vspace{2mm}

Address (both authors): \small

University 'Al.I.Cuza' of Ia\c si\\
Faculty of Mathematics\\
Bd. Carol I, no.11\\
700506 Ia\c si\\
Romania\\
e-mail (M.I.Munteanu): marian.ioan.munteanu@gmail.com\\
e-mail (A.I.Nistor): ana.irina.nistor@gmail.com \hspace{40mm}
Received,

\normalsize

\end{document}